\documentclass[12pt]{amsart} 
\usepackage{mathrsfs}
\usepackage{amssymb,amsmath,amsthm,color}
\usepackage{graphicx}
\usepackage{hyperref}
\usepackage{url} 
\usepackage{setspace}
\usepackage{mathtools}
\usepackage[inline]{enumitem}  
\usepackage[margin=1in]{geometry}
\usepackage{comment}
\usepackage{cite}
\numberwithin{equation}{section}
\newtheorem{theorem}{Theorem}[section]

\newtheorem{Definition}{Definition}[section]
\newtheorem{Remark}{Remark}[section]
\newtheorem{proposition}{Proposition}[section]

\newcommand{\R}{\mathbb R}

\begin{document}
\baselineskip16pt
\title[The  Klein-Gordon-Hartree equation in modulation spaces]{The Global well-posedness for Klein-Gordon-Hartree equation in modulation spaces}
\author{Divyang G. Bhimani}
\address{Department of Mathematics, Indian Institute of Science Education and Research, Dr. Homi Bhabha Road, Pune 411008, India}
\email{divyang.bhimani@iiserpune.ac.in}

\date{}
\maketitle
\begin{abstract}Modulation  spaces have received considerable interest recently as it is the natural function spaces to consider  low regularity Cauchy data for several nonlinear evolution equations.   We  establish global well-posedness for  3D Klein-Gordon-Hartree equation  $$u_{tt}-\Delta u+u + ( |\cdot|^{-\gamma} \ast |u|^2)u=0$$
 with initial  data in  modulation spaces $M^{p, p'}_1 \times M^{p,p}$ for $p\in \left(2, \frac{54	}{27-2\gamma} \right),$  $2<\gamma<3.$  We implement  Bourgain's  high-low  frequency decomposition method to  establish global well-posedness,   which was earlier used for classical  Klein-Gordon equation.
 This is the first  result on low regularity for Klein-Gordon-Hartree equation with large initial data in modulation spaces (which do not coincide with Sobolev spaces).
 \end{abstract}
%\tableofcontents
\section{Introduction}
 We study the Cauchy problem for Klein-Gordon-Hartree equation
 \begin{equation}\label{kgh}
 \begin{cases} u_{tt}-\Delta u+u + (V \ast |u|^2)u=0\\
 u|_{t=0}=f, \quad \partial_tu|_{t=0}=g
 \end{cases} (t,x) \in \mathbb R \times \R^3.
 \end{equation}
 Here, $u(t,x)\in \mathbb C$ and $V$ is a Hartree kernel given by 
 \[ V(x)= \frac{1}{|x|^{\gamma}} \quad (0<\gamma<3).\]
 
 The equation \eqref{kgh} has  appeared in the context of various physical models,  e.g.  one-component plasma,   relativistic  particles,   quantum semiconductor devices,  etc.. See \cite{JeanS}  and the references therein for a general review.  
 The case  $\gamma=1$ in $V$  is known as  Coulomb  potential,  which play an important  role in quantum mechanics.   Compared with the classical Klein-Gordon equation with local  non-linearity $|u|^{p-1}u$,  the  nonlinearity 
 $ (V\ast |u|^2)u$ is non local.

Formally,  the solution of \eqref{kgh} conserves the energy
\begin{eqnarray*}
E(u(t),  \partial_t u(t)) &  = & \frac{1}{2} \int_{\R^3} \left( |\partial_t u(t,x)|^2 + |\nabla u (t,x)|^2 + |u(t,x)|^2 \right) dx\\
& \   \ & + \frac{1}{4}\int_{\R^3} \int_{\R^3} \frac{|u(t,x)|^2 |u(t,y)|^2}{|x-y|^{\gamma}} dxdy\\
& = &  E (f, g)
\end{eqnarray*}
and the momentum 
\begin{eqnarray}
P(u)(t)= \int_{\R^3} u_t(t, x) \nabla u(t,x) dx = P(u)(0).
\end{eqnarray}

In the early 1980s,  G.  Perla  and W. Strauss \cite{Perla1982}   have initiated 
the study of  \eqref{kgh}.  Later W.  W. Strauss \cite[Theorem 19]{s2} (cf \cite[Theorem 11]{s1}) proved global well-posedness and scattering  for \eqref{kgh} with  $V\in L^r \ (1\leq r \leq 3/2)$ in  the energy space  $H^1\times L^2.$  Mochizuki\cite{mochizuki1989JMSJ}
  proved  small data  global well-posedness and scattering  for \eqref{kgh} with $2\leq \gamma <3$  in   $H^1\times L^2.$  Later Miao-Zheng \cite[Theorem 1.1]{Miao2014}  generalized this  result  for the   general initial data in $H^1\times L^2$ for subcritical non-linearity (i.e.  $2<\gamma <3$).   In \cite{QMiao2015}, Miao-Zheng  established global well-posedness for defocusing \eqref{kgh} in higher dimensions  in $H^1\times L^2.$
  While  Miao and Zhang  \cite[Theorem 1.1]{Miao2011} proved  global well-posedness for \eqref{kgh} with low regularity Cauchy data, i.e.  data in  below energy space $H^{s} \times H^{s-1} \ (\gamma/4 <s<1)$ for $2<\gamma<3.$  See also \cite{Miao2014}.   We note that these results did not cover Coulomb  potential.  We also remark that  all previous authors have studied  \eqref{kgh} on $L^2-$based Sobolev spaces.  Mainly because   Klein-Gordon  semigroup $e^{it(I-\Delta)^{1/2}}$   fails to be bounded  on $L^p$ if $p\neq2.$ Hence we cannot expect to solve  equations \eqref{kgh} in  $L^p-$spaces for $p\neq 2$.   The question arises if it is possible to remove $L^2$ constraint and consider   \eqref{kgh}  in function spaces which are not  just $L^2$ based.

In this paper we study \eqref{kgh} with Cauchy data in the modulation spaces  $M^{p,q}_s-$   which are  based on  mixed    $\ell_s^q(L^p)$norms and 
contain rougher  Cauchy data  than  in any given fractional Bessel potential space (see e.g.  \cite[Chapter 6]{KassoBook}).   In order to define these spaces, we briefly  introduce some notations.
Let   $\rho \in \mathcal{S}(\mathbb R^d)$ (Schwartz space),  $\rho: \mathbb R^d \to [0,1]$  be  a smooth function satisfying   $\rho(\xi)= 1 \  \text{if} \ \ |\xi|_{\infty}\leq \frac{1}{2} $ and $\rho(\xi)=
0 \  \text{if} \ \ |\xi|_{\infty}\geq  1.$ Here,  $|\xi|_{\infty}=\max\{ | \xi_i | : \xi= (\xi_1,..., \xi_d)\}.$ Let  $\rho_k$ be a translation of $\rho,$ that is,
$ \rho_k(\xi)= \rho(\xi -k) \ (k \in \mathbb Z^d).$
Denote 
$\sigma_{k}(\xi)= \frac{\rho_{k}(\xi)}{\sum_{l\in\mathbb Z^{d}}\rho_{l}(\xi)}\quad (k \in \mathbb Z^d).$
The frequency-uniform decomposition operators can be  defined by 
$$\square_k = \mathcal{F}^{-1} \sigma_k \mathcal{F}$$
where $\mathcal{F}$ and $\mathcal{F}^{-1}$ denote the Fourier and inverse Fourier transform respectively.   The weighted modulation spaces  $M^{p,q}_s(\R^d) \ (1 \leq p,q \leq \infty, s \in \R)$ is defined as follows:
 \[ M^{p,q}_s(\R^d)= \left\{ f \in \mathcal{S}'(\R^d): \|f\|_{M^{p,q}_s}:   \left\| \|\square_kf\|^q_{L^p_x} (1+|k|)^{s} \right\|_{\ell^q_k}< \infty  \right\} . \]
See also Remark \ref{edm}.  For $s=0,$ we write $M^{p,q}_0(\R^d)= M^{p,q}(\R^d).$ For $p\in [1, \infty]$,  we denote $p'$ the Hilbert conjugate, i.e.  $\frac{1}{p}+\frac{1}{p'}=1.$
It is known that 
\[M^{2,2}_s(\R^d)= H^s(\R^d) \ \text{(Sobolev space)} \quad \text{and} \quad  M^{p,p'}(\R^d) \hookrightarrow  L^p(\R^d)  \ \text{for}  \   p \in [2,  \infty].\]
In the last two decades $M^{p,q}_s-$spaces have   been  considered extensively  as low regularity  Cauchy data spaces  for several  nonlinear evolution  equations, see  \cite{BhimaniHartree-Fock, Bhimani2016,  Kasso2009, BhimaniNorm, Wang2006, Wang2007, WangBook, BhimaniEjDE, KassoBook,  LeonidIn, Pattakos2021}.  We shall briefly mention some of them here.
In the seminal work of    
J. Ginibre and G. Velo \cite[Proposition 3.2]{GVMZ},   they proved global well-posedness in  $H^1\times L^2$ for  the classical  Klein-Gordon  (CKG) equation,  i.e.  \eqref{kgh}  with local nonlinearity $|u|^2u$ (in other words,  \eqref{kgh} with $V$ is the  Dirac delta distribution).   While in the seminal work of B. Wang and H. Hudzik  \cite[Theorem 1.6]{Wang2007},    they established  low regularity  small data global well-posedness for the CKG equation in some  modulation spaces $M^{2,1}_{s} \times M^{2,1}_{s-1}.$ See also \cite{GZ2014}.  In \cite[Theorem 1.3]{Kasso2009},   B\'enyi and Okoudjou proved  local well-posedness for  CKG
 in $M^{p,1}_s.$  See also \cite[Theorem 4.8]{Bhimani2016}.   Recently L. Chaichenets and N.  Pattakos  \cite{Pattakos2021} proved global well-posedness for CKG equation in $M_1^{p, p'}\times M^{p, p}$ for $p$ sufficiently close to $2.$

 Now we state our main results.

 \begin{theorem}[Local well-posedness]\label{mt1} Let $(f, g)\in H^1(\R^3)\times L^2(\R^3)+M_1^{p_{\gamma}, p'_{\gamma}}(\R^3)\times M^{p_{\gamma},  p'_{\gamma}}(\R^3),$ where $p_{\gamma}= \frac{18}{9-2\gamma}$ and $0<\gamma<3.$ Then,  there exists $T=T(\| (f, g)\|_{H^1\times L^2+M_1^{p_\gamma,p'_\gamma}\times M^{p_\gamma, p'_\gamma}})>0$ such that \eqref{kgh} has unique  solution $$u\in C([0, T],  H^1(\R^3))+C([0, T], M_1^{p_\gamma,p'_{\gamma}}(\R^3))$$ with $\partial_tu\in u\in C([0, T],  L^2(\R^3))+C([0, T], M^{p_\gamma, p'_\gamma}(\R^3)).$ Moreover,  if $T^*$ is the maximum time of existence  then the blowup alternative holds,  i.e.,   if $T^*<\infty,$ then
 \[ \limsup_{t\to T^*} \| \left( u(t, \cdot), \partial_t u(t, \cdot) \right)\|_{H^1\times L^2+M_1^{p_\gamma,p'_\gamma}\times M^{p_\gamma,  p'_\gamma}}= \infty. \]
 \end{theorem}
For $\gamma \in (0,3)$ we have $p_\gamma \in (2,6).$
Notice that $M^{p_\gamma,  p'_\gamma}(\mathbb R^3) \hookrightarrow L^{p_\gamma}(\mathbb R^3)$  is a sharp embedding and up to now we cannot get the local well-posedness of \eqref{kgh} in $L^{p_\gamma}(\R^3)$ but in $M^{p_\gamma,  p'_\gamma}(\mathbb R^3).$

\begin{Remark} We have several comments for  Theorems \ref{mt1}.
\begin{enumerate}
 \item The proof of  local well-posedness  is based on  Banach contraction principle.  Specifically,   we show that the integral version of \eqref{kgh},  i.e.  the operator $\Phi$ in \eqref{fpa},  is a contraction on some suitable ball. 
 \item  In order to  control the linear  part  in \eqref{fpa},  we use the fact that Klein-Gordon semigroup is bounded on modulation spaces (see  Proposition \ref{bm} below).   While in order to control the nonlinear part in \eqref{fpa}, we make use of  Strichartz estimates,   Hardy-Littlewood Sobolev  inequality and   Sobolev  embedding.
 \item We shall see that  the life span of the local solution  depends on the norm of the initial data,  and  so one can follow standard procedure to establish blowup alternative.
 \item In order to treat  Hartree non-linearity,  we consider  $p_{\gamma}= \frac{18}{9-2\gamma},$which turned out to be compatible with the several inequalities we used.  See \eqref{air}.
\end{enumerate}
\end{Remark}

\begin{theorem}[Global well-posedness]\label{mt2}
Let real initial data $(f,g)\in M_1^{p,p'}(\R^3)\times M^{p, p'}(\R^3)$  where $p\in \left(2, \frac{54}{27-2\gamma} \right)$ with $2<\gamma<3.$
Then the local solution constructed in Theorem \ref{mt1} of \eqref{kgh}  extends globally and lies in
\[u \in C(\R, H^1(\R^3))+ C(\R, M_1^{p_\gamma,  p'_\gamma}(\R^3)) \]
with $\partial_tu\in C(\R,  L^2(\R^3))+ C(\R, M^{p_\gamma,  p'_\gamma }(\R^3)).$
\end{theorem}  
Theorem \ref{mt2} is new and complement known results for \eqref{kgh}  in \cite{Perla1982, Perla1983, Miao2011, Miao2014,  QMiao2015, BhimaniEjDE}.  
Exploiting  the decay estimate (Proposition \ref{wp}),   Bhimani \cite[Theorem 1.1]{BhimaniEjDE}  proved  small data global well-posedness for \eqref{kgh} for in some modulation space.     While Theorem \ref{mt2}  establishes global well-posedness for any data (large) in   some modulation spaces.

\begin{Remark} We  briefly discuss proof strategy and imposed hypothesis of Theorem \ref{mt2}.
\begin{enumerate}
\item The proof of global well-posedness is inspired from  Bourgain's  high-low frequency decomposition method (see \cite[Section IV.2]{Bourgain1999}) and  Chaichenets-Pattakos's  Theorem 1.2 in  \cite{Pattakos2021})-where they considered CKG equation.  See also \cite[Section 4.2]{TzirakisBook} and \cite[Section 3.9]{TaoBook}. 
\item  Using interpolation theory decompose initial data into two parts:$$(f,g)=(f_N, g_N)+ (f^N, g^N)\in X_0 + Y_0,$$ as in \eqref{spbt},  say $(f_N, g_N)$ good part (low frequencies) and $(f^N, g^N)$ bad part (high frequency).

\item In view of  blow-up alternative,   it is sufficient  to prove that $$ \| \left( u(t, \cdot), \partial_t u(t, \cdot) \right)\|_{H^1\times L^2+M_1^{p_\gamma,p'_\gamma}\times M^{p_\gamma,  p'_\gamma}}< \infty$$  on \textit{any} bounded time interval $I=[0, T]$ to ensure the  global existence.

\item To this end,   we  convert \eqref{kgh} into Schr\"odinger equation $$iv_t-Bv-B^{-1} \left( [V \ast | \text{Re} \ v|^2 ] \text{Re} \ v \right) =0,  \quad v(0) = F_N + F^N \in H^1 + M^{p_\gamma,  p'_\gamma }$$ via change of variables $B= (Id-\Delta)^{1/2}.$ See \eqref{bit}- \eqref{mg}.  Since the linear part is bounded on modulation space (Propositions \ref{fp} and \ref{bm}),   we are left to control the energy space  $H^1-$norm of  $\tilde{v}(t)=(v(t)- e^{-itB}F^N)$.   This we shall do this using the Hamiltonian $H$.  See \eqref{hamil}. 
\item  Define $I(t)=H(\tilde{v}(t)).$ We shall show that $I(0)\lesssim N^{\frac{4\theta}{1-\theta}}.$  The restriction on $2<\gamma$  comes at this stage.  See \eqref{s1}.   We push this forward to $[0,T]$ by looking at $I(t)$.   

\item  We shall observe that $T\to \infty$ as $N\to \infty;$ which will yield global existence.  This will impose conditions on $p.$  In fact,  the condition on $p$ comes due to \eqref{rc}- where we make use of the identity \eqref{ipu} in order to obtain the upper bound for $p$ in terms of  Hartree exponent $\gamma.$
\end{enumerate}
\end{Remark}

\begin{Remark}  Our method of proof  should be applied to treat  higher dimension $d>3$ and also to more general Choquard  type nonlinearity  $(|\cdot|^{-\gamma}\ast |u|^{\alpha})|u|^{\alpha-2}u$ and its variants.  See \cite{JeanS}.   In the present paper  we have confined ourselves to the size of the physical space $d = 3$ and $\alpha=2$ for the presentation of clear ideas.
\end{Remark}

\begin{Remark} \label{edm} 

In the early 1980s  Feichtinger  \cite{Feih83} introduced a  class of Banach spaces,  the so called modulation spaces, which allow a measurement of space variable and Fourier transform variable of a function or distribution on $\mathbb R^d$ simultaneously using the short-time Fourier transform (STFT).  The  STFT  of a tempered distribution $f \in \mathcal{S}'(\R^d)$ with respect to a window function $0\neq g \in {\mathcal S}(\R^d)$ (Schwartz space) is defined by
\begin{eqnarray*}\label{stft}
V_{g}f(x,y)= \int_{\mathbb R^{d}} f(t) \overline{g(t-x)} e^{-2\pi i y\cdot t}dt,  \  (x, y) \in \mathbb R^{2d}
\end{eqnarray*}
 whenever the integral exists.  It is known \cite{Feih83},  \cite[Proposition 2.1]{Wang2007} that 
\begin{eqnarray*}
\|f\|_{M^{p,q}_s}\asymp   \left\| \|V_gf(x,y)\|_{L^p_x} (1+ |y|^2)^{s/2} \right\|_{L_y^q}.
\end{eqnarray*}
\end{Remark}
\section{Preliminaries}\label{mp}
\begin{proposition}[\cite{KassoBook, WangBook},  see Theorem 9.9 in  \cite{Brez}]\label{up} Let $p,q \in [1,  \infty]$ and $s\in \R.$ 
\begin{enumerate}
\item $\mathcal{S}(\R^d) \hookrightarrow M^{p,q}(\R^d)\hookrightarrow S'(\R^d).$
\item \label{eml} $M^{p,p'}(\R^d) \hookrightarrow L^p(\R^d)$ for $p \in [2, \infty].$
\item \label{up1} $H^1(\R^d)\hookrightarrow L^p(\R^d)$ for all $p\in \left[ 2, \frac{2d}{d-2}\right]$  (Sobolev embedding).
\end{enumerate}
\end{proposition}
Recall that equation \eqref{kgh} have the following
equivalent form
\[u(t, x)= K'(t)f(x) + K(t)g(x) -\mathcal{B}f(u),\]
where 
\[K(t)=\frac{\sin t (I-\Delta)^{1/2}}{ (I-\Delta)^{1/2}}, \  \  K'(t)= \cos t (I-\Delta)^{1/2}, \ \  \mathcal{B}= \int_0^t K(t-\tau) \cdot d\tau.\]
Notice that the (semi-) group given by 
$$\mathbb K (t) =\left(\begin{matrix}
  K'(t) & K(t)\\
  (\Delta-I)K(t) & K'(t)
\end{matrix}\right) \quad  \forall t \in \R$$
is isometric on $H^1(\R^d) \times L^2 (\R^d)$ with domain of its generator being $\mathcal{D}=H^2(\R^d) \times H^1(\R^d)$.
\begin{proposition}[See Proposition 4.2 in  \cite{Wang2007}]\label{wp} Let $G(t)= e^{it (I-\Delta)^{1/2}} \ (t\in \mathbb R).$
\begin{enumerate}
\item Let $s\in \mathbb R, 2\leq p \leq \infty,$ $1 \leq  q < \infty, \theta \in [0,1],$ and
$2 \sigma (p)= (d+2) \left( \frac{1}{2} - \frac{1}{p}\right).$
Then we have 
\[  \|G(t) f\|_{M^{p,q}_s} \lesssim (1+ |t|)^{-d\theta  (1/2-1/p)} \|f\|_{M^{p',q}_{s+ \theta 2 \sigma(p)}}. \]
\item  \label{kb} Let  $s\in \mathbb R $ and $ 1\leq p, q  \leq \infty.$   Then we have 
\[  \|G(t)f\|_{M^{p,q}_s} \leq C (1+ |t|)^{d \left| 1/2 -1/p \right|}  \|f\|_{M^{p,q}_s}.\]
\end{enumerate}
\end{proposition}

\begin{proposition}[See Theorem 6.6 in \cite{Feih83},  Proposition 2.4 in \cite{Wang2007}] \label{fp} Let $p,q \in [1, \infty]$ and $s, \sigma \in \R.$ Then the map $(Id-\Delta)^{\frac{\sigma}{2}}$ is an isomorphism from $M^{p,q}_s(\R^d)$ onto $M^{p,q}_{s-\sigma}(\R^d)$.
\end{proposition}
\begin{proposition}[Uniform boundedness of Klein-Gordon propagator, See Corollary 3.2 in \cite{EFwave}]\label{bm} Let $p,q \in [1, \infty]$ and $s\in \mathbb R$. Then for every $T>0$ there  is a constant $C_T=C_T(p,q,s)$ such that 
\[\|K(t)f\|_{M^{p,q}_{s+1}}+\|K'(t)f\|_{M^{p,q}_s}\leq C_T \|f\|_{M^{p,q}_s} \]
for all $f\in M^{p,q}_s(\mathbb R^d)$ and all $|t|\leq T.$
\end{proposition}
\begin{Definition}\label{kgadp} We say that  the pair $(q,r) \in [2, \infty] \times [2, \infty)$ is (wave-)admissible if 
\begin{eqnarray}
\frac{1}{q} + \frac{d-1}{2r} \leq \frac{d-1}{4}.
\end{eqnarray}
We will denote by $q_a(r)$ the unique solution of the gap condition
\[ \frac{1}{q} + \frac{d}{r}= \frac{d}{2}-1. \]
\end{Definition}
In Theorem \ref{st} we recall well known  Strichartz estimates (see 
\cite[Corollary 1.3]{KeelTao1998},   \cite[Estimate 1.1]{Zheng2019},
and \cite[Section 2.2]{GrillakisBook} for details). We use  $L_t^{r} ([0,T], X)$ to denote the space time norm
\[\|u\|_{L^r_TX}=\|u\|_{L^{r}_t([0,T] X)}=  \left( \int_{0}^T \|u\|^r_{X}  dt \right)^{1/r}.\]
where  $X$ is a Banach space.  

\begin{theorem}[Strichartz estimates]\label{st} If $u$ is a solution to \eqref{kgh} with initial data $(f,g) \in H^1(\R^d) \times L^2(\R^d),$ then we have 
\begin{eqnarray*}
\|u\|_{L^q_TL^r}+\|u\|_{C_TH^1}+\|\partial_tu\|_{C_TL^2}\lesssim \|f\|_{H^1}+\|g\|_{L^2}+\|(V \ast |u|^2)u\|_{L_T^1L^2}
\end{eqnarray*}
for every admissible pair $(q,r)$ that satisfies the  gap condition in Definition \ref{kgadp}, i.e. 
\[q=\left( d \left( \frac{1}{2}- \frac{1}{r} \right) -1 \right)^{-1},  \quad r \in \left[  \left( \frac{1}{2}- \frac{1}{d-1} \right)^{-1},  \frac{2(d+1)}{d-3} \right].  \]
\end{theorem}

We briefly recall some useful fact in interpolation theory.  We refer to \cite{Interpolation} for a comprehensive introduction.  We say that the Banach spaces $X, Y$ form the interpolation couple $\{X, Y\}$ if there is a topological  Hausdroff space $\mathcal{V}$ such that $X, Y \subset \mathcal{V}$ and $X, Y \hookrightarrow \mathcal{V}.$  For an interpolation couple $\{ X, Y \}, t \in (0, \infty)$ and $z\in X+Y,$ we set  (the $K-$functional)
\[ K(t,z)= \inf_{{z=x+y} \atop {x\in X, y \in Y}} \left( \|x\|_{X} +t \|y\|_{Y} \right). \]
\begin{Definition} Let $\{X, Y\}$ be an interpolation couple  and $\theta \in (0, 1).$ Put $\|z\|_{(X, Y)_{\theta, \infty}}= \sup_{t>0} \left[ t^{-\theta} K(t,z)\right].$The real interpolation space is defined by
\[ (X, Y)_{\theta, \infty} = \left\{ z \in X+Y: \|z\|_{(X,Y)_{\theta, \infty}} < \infty \right\}.\]
\end{Definition}
\begin{proposition}[See Lemma 16 in  \cite{Pattakos2021}]\label{ipp}   Let $\{X, Y\}$ be an interpolation couple, $\theta \in (0, 1)$ and $f\in [X, Y]_{\theta}.$ Then for any $N\in \R^+$ there is an $f_N \in X$ and $f^{N}\in Y$ such that
\[ \|f_N\|_{X} \lesssim N^{\frac{\theta}{1-\theta}}\|f\|_{[X, Y]_{\theta}} \quad \text{and} \quad  \|f^{N}\|_{Y} \lesssim  \frac{1}{N} \|f\|_{[X, Y]_{\theta}}.\]
\end{proposition}
\begin{proposition}[Hardy-Littlewood-Sobolev inequality] \label{hls}Assume that  $0<\gamma< d$ and $1<p<q< \infty$ with
$\frac{1}{p}+\frac{\gamma}{d}-1= \frac{1}{q}.$
Then 
$\||\cdot|^{-\gamma}\ast f\|_{L^q} \leq C_{d,\gamma, p} \|f\|_{L^p}.$
\end{proposition}

\begin{proposition}[Complex interpolation,  see Theorem 6.1 D in \cite{Feih83}] \label{ci} Let $p_0, p_1, q_0, q_1 \in [1, \infty]$  such that $(q_0, q_1) \neq (\infty,  \infty).$ Furthermore,  let $s_0, s_1 \in \R$ and $\theta \in (0,1).$  Define $s= (1-\theta)s_1 + \theta s_2 \in \R$ and $p \in [1, \infty]$ via
\[ \frac{1}{p}= \frac{1-\theta}{p_0} + \frac{\theta}{p_1}. \]
Finally,  define $q\in [1,  \infty)$ via 
\[ \frac{1}{q}= \frac{1-\theta}{q_0} + \frac{\theta}{q_1}. \]
Then
\begin{equation*}
\left[M^{p_0,  q_0}_{s_0}(\R^d),  M^{p_1, q_1}_{s_1}(\R^d)\right]_{\theta}= M^{p,q}_s(\R^d).
\end{equation*}
\end{proposition}

\section{Proof of Theorems \ref{mt1} and \ref{mt2}}
\begin{proof}[\textbf{Proof of Theorem \ref{mt1}}]
We may decompose initial data
 $$(f, g)\in X:= H^1(\R^3)\times L^2(\R^3)+M_1^{p_{\gamma}, p'_{\gamma}}(\R^3)\times M^{p_{\gamma},  p'_{\gamma}}(\R^3)$$ as follows
\begin{eqnarray}\label{spdt}
(f,g)=(f_0,g_0)+(\tilde{f_0}, \tilde{g_0})= (f_0+\tilde{f_0}, g_0+ \tilde{g_0})
\end{eqnarray} 
where
\[(f_0,g_0)\in H^1(\R^3)\times L^2(\R^3) \quad \text{and} \quad (\tilde{f_0}, \tilde{g_0}) \in M_1^{p_{\gamma}, p'_{\gamma}}(\R^3)\times M^{p_{\gamma},  p'_{\gamma}}(\R^3).\]
We present  an argument for the first coordinate for the solution $(u, \partial_tu)$, i.e. for $u$  only.   As taking  Theorem \ref{st} into account,  the argument for $\partial_tu$ is similar. 
Consider the Banach space $$X(T)=X_1(T)+X_2(T)$$ where
\[X_1(T)=L^{\infty}([0, T], H^1(\R^3)) \quad \text{and} \quad X_2(T)=L^{\infty}([0,T], M^{p_{\gamma},  p'_{\gamma}}(\R^3)).\]
The norm in $X(T)$ is given by 
\[ \|u\|_{X(T)}= \inf_{{u=u_1+u_2}\atop{u_1\in X_1(T), u_2\in X_2(T)}} \left(\|u_1\|_{X_1(T)} + \|u_2\|_{X_2(T)} \right).\]
We write \eqref{kgh} as
\begin{eqnarray}\label{fpa}
u(t)= K'(t)f + K(t)g - \int_0^t K(t-\tau) (V\ast |u|^2)u(\tau) d\tau:=\Phi(u)(t).
\end{eqnarray}

We consider the ball $B(R,T)=\{ u\in X(T): \|u\|_{X(T)} \leq R \}$ for some $R$ and $T>0$ (to be chosen below).
We  shall show that for the operator $\Phi$ is a contraction on $B(R,T)$ for some $R$ and $T>0.$  We shall first show that $\Phi$ is a self map on $B(R,T).$ We  fix $u\in B(R, T)$ and consider decomposition $u=v+w,$ $v\in X_1(T)$ and $w\in X_2(T).$

Without loss of generality we may assume that $T\leq 1.$
By \eqref{spdt},  Theorem \ref{st} and Proposition \ref{bm},  we obtain 
\begin{eqnarray*}
&& \|K'(t)f + K(t)g\|_{X(T)}\\
 &= & \| K'(t)f_0 + K(t)g_0+K'(t)\tilde{f_0} + K(t)\tilde{g_0}\|_{X(T)}\\
  & \lesssim  & \|K'(t)f_0 + K(t)g_0 \|_{X_1(T)}+ \|K'(t)\tilde{f_0} + K(t)\tilde{g_0}\|_{X_2(T)}\\
 & \lesssim & \|f_0\|_{H^1}+\|g_0\|_{L^2}+ C(T)\left( \|\tilde{f_0}\|_{M_1^{p_{\gamma},  p'_{\gamma}}} +  \|\tilde{g_0}\|_{M^{p_{\gamma},  p'_{\gamma}}}   \right)\\
 & \lesssim & \|f_0\|_{H^1}+\|g_0\|_{L^2}+  \|\tilde{f_0}\|_{M_1^{p_{\gamma},  p'_{\gamma}}} +  \|\tilde{g_0}\|_{M^{p_{\gamma},  p'_{\gamma}}}.
\end{eqnarray*}
Since the decomposition of a given data is arbitrary (\eqref{spdt}),   we obtain
\begin{eqnarray}\label{lec}
\|K'(t)f + K(t)g\|_{X(T)} \lesssim \| (f, g)\|_{X}.
\end{eqnarray}
We now choose $R\approx 2 \|(f,g)\|_{X}.$ 
By Proposition \ref{up},  we have 
\begin{eqnarray}\label{ce}
 X(T) \hookrightarrow L^{\infty}([0,T], L^{p_{\gamma}} (\R^3)).
\end{eqnarray}
Notice that 
\begin{eqnarray}\label{uie}
 \frac{\gamma}{9}+\frac{9-2\gamma}{18}=\frac{1}{2},  \quad  \frac{2}{3}+\frac{1}{3}=1\quad \text{and} \quad \frac{9-2\gamma}{9}+\frac{\gamma}{3}-1=\frac{\gamma}{9}.
\end{eqnarray}
By  Theorem \ref{st},  H\"older  inequality,   Proposition \ref{hls} (with $p= \frac{9}{9-2\gamma},  q=\frac{9}{\gamma}, d=3$)  and \eqref{ce},  we obtain 
\begin{eqnarray}\label{air}
\left\|  \int_0^t K(t-\tau) (V\ast |u|^2)u(\tau) d\tau \right\|_{L^{\infty}_TH^1} & \lesssim & \|(V\ast |u|^2) u\|_{L^1_TL^2}= \left\| \|(V\ast |u|^2) u\|_{L^2} \right\|_{L^{1}_T} \nonumber \\
& \lesssim &   \left\| \|(V\ast |u|^2) \|_{L^{\frac{9}{\gamma}}} \|u \|_{L^{\frac{18}{9-2\gamma}}} \right\|_{L^{1}_T} \nonumber \\
& \lesssim &  \|V\ast |u|^2\|_{L^{\frac{3}{2}}_TL^{\frac{9}{\gamma}}}\|u\|_{L^3_TL^{\frac{18}{9-2\gamma}}} \nonumber  \\
& \lesssim &  \left\| \| u\|^2_{L^\frac{18}{9-2\gamma}} \right \|_{L^{\frac{3}{2}}_T} \|u\|_{L^3_TL^{\frac{18}{9-2\gamma}}}  =  \|u\|^3_{L^3_TL^{\frac{18}{9-2\gamma}}}  \nonumber \\
& \lesssim & T \|u\|_{L_T^{\infty}L^{\frac{18}{9-2\gamma}}}^3=T \|u\|_{L_T^{\infty}L^{p_{\gamma}}}^3\\
& \lesssim & T \|u\|_{X(T)}^3 \nonumber \\
& \lesssim & TR^3. \nonumber 
\end{eqnarray}
Using this together with \eqref{lec} and taking $T=\frac{1}{2R^2}=\frac{1}{8\|(f,g))\|_{X}^2}$,  we conclude that   $\Phi$ is self map on $B(R,T).$ Next,  we shall prove $\Phi$ is a contraction on $B(R, T).$
Let $u_1, u_2 \in X(T).$  Recall the identity
\begin{eqnarray}\label{pehn}
(V\ast |u_1|^{2})u_1- (V\ast |u_2|^{2})u_2= (V\ast |u_1|^{2})(u_1-u_2) + (V \ast (|u_1|^{2}- |u_2|^{2}))u_2.
\end{eqnarray}  
Using \eqref{pehn} and performing similar arguments  as above,  we obtain
\begin{eqnarray*}
&& \|\Phi(u_1)-\Phi(u_2)\|_{L^{\infty}_TH^1} \\
& \lesssim & \|(V\ast |u_1|^{2})(u_1-u_2)\|_{L^1_TL^2}+ \|(V \ast (|u_1|^{2}- |u_2|^{2}))u_2\|_{L^1_TL^2}\\
& \lesssim & T \|u_1\|_{L^{\infty}_TL^{p_{\gamma}}}^2\|u_1-u_2\|_{L^{\infty}_TL^{p_\gamma}} + \left\|  \||u_1|^2-|u_2|^2\|_{L^{\frac{9}{9-2\gamma}}} \right\|_{L_T^{\frac{3}{2}}} \|u_2\|_{L^3_TL^{\frac{18}{9-2\gamma}}}.
\end{eqnarray*}
By Cauchy–Schwartz  and H\"older inequalities,  we obtain
\begin{eqnarray*}
 && \left\|  \||u_1|^2-|u_2|^2\|_{L^{\frac{9}{9-2\gamma}}} \right\|_{L_T^{\frac{3}{2}}}\\
  & \lesssim & \left\| \|(u_1-u_2)\bar{u_1}\|_{L^{\frac{9}{9-2\gamma}}}+ \|(\bar{u_1}-\bar{u_2})u_2\|_{L^{\frac{9}{9-2\gamma}}}\right\|_{L_T^{\frac{3}{2}}}\\
 & \lesssim & \left\| \|\bar{u_1}\|_{L^{\frac{18}{9-2\gamma}}}\|u_1-u_2\|_{L^{\frac{18}{9-2\gamma}}}+ \|(\bar{u_1}-\bar{u_2})\|_{L^{\frac{18}{9-2\gamma}}} \|u_2\|_{L^{\frac{18}{9-2\gamma}}} \right\|_{L_T^{\frac{3}{2}}}\\
 & \lesssim & \|\bar{u_1}\|_{L^3_TL^{\frac{18}{9-2\gamma}}}\|u_1-u_2\|_{L^3_TL^{\frac{18}{9-2\gamma}}}+ \|\bar{u_1}-\bar{u_2}\|_{L^3_TL^{\frac{18}{9-2\gamma}}} \|u_2\|_{L^3_TL^{\frac{18}{9-2\gamma}}}. 
\end{eqnarray*}
Thus, we  have 
\begin{align*}
& \|\Phi(u_1)-\Phi(u_2)\|_{X_1(T)}\\
 & \lesssim   T \left( \|u_1\|_{L^{\infty}_T L^{p_{\gamma}}}^2+\|u_1\|_{L^{\infty}_T L^{p_\gamma}}\|u_2\|_{L^{\infty}_T H^1}+\|u_2\|^2_{L^{\infty}_T L^{p_\gamma}} \right) \|u_1-u_2\|_{L^{\infty}_T L^{p_\gamma}}\\
 & \lesssim  TR^2 \|u_1-u_2\|_{X(T)}
\end{align*}
By choosing  $T$ further small enough (if necessary),   it follows that $\Phi$ is a contraction on $B(R,T)$. This completes the proof of Theorem \ref{mt1}.
\end{proof}
\begin{proof}[\textbf{Proof of Theorem \ref{mt2}}]
Let $(f,g)\in M_1^{p,p'}(\R^3)\times M^{p, p'}(\R^3).$   By   Proposition \ref{ci}, we have 
\begin{eqnarray}
M^{p,p'}_1(\R^3)=\left[ H^1(\R^3), M^{p_\gamma, p'_\gamma}(\R^3) \right]_{\theta}, \quad M^{p, p'}(\R^3)= \left[ L^2(\R^3), M^{p_\gamma, p'_\gamma}(\R^3) \right]_{\theta}
\end{eqnarray}
where $\theta\in (0,1)$ is fixed by 
\begin{eqnarray}\label{ipu}
\frac{1}{p}= \frac{1-\theta}{2}+ \frac{\theta}{p_\gamma}.
\end{eqnarray}
By Proposition \ref{ipp}, there is a decomposition 
\begin{eqnarray} \label{spbt}
(f,g)=(f_N, g_N)+ (f^N, g^N)\in H^1(\R^3)\times L^2(\R^3)+M_1^{p_\gamma,p'_\gamma}(\R^3)\times M^{p_\gamma, p'_\gamma}(\R^3)
\end{eqnarray}
with
\begin{eqnarray}\label{tbu}
\|f_N\|_{H^1}, \|g_N\|_{L^2} \lesssim N^{\frac{\theta}{1-\theta}},  \quad  \|f^N\|_{M^{p_\gamma, p'_\gamma}_1}, \|g^N\|_{M^{p_\gamma, p'_\gamma}} \lesssim \frac{1}{N}
\end{eqnarray}
for any $N>0$ (to be fixed later, depending on time $T$).

Note that in view of Theorem \ref{mt1} for real initial data $(f,g)$ there is a real-valued local solution $u$ to \eqref{kgh}.
Now we shall see that the solution constructed before in Theorem \ref{mt1} is global in time.  In fact, in view of
blowup alternative,  to prove Theorem \ref{mt2},  it suffices to prove that the norm
\[ \| \left( u(t, \cdot), \partial_t u(t, \cdot) \right)\|_{H^1\times L^2+M_1^{p_\gamma,p'_\gamma}\times M^{p_\gamma,  p'_\gamma}} \]
remain finite  on \textit{any} bounded time intervals $I=[0,T]$.  We may now perform change of variable as in \cite[Equation IV. 2.4]{Bourgain1999}. 
Specifically,  we let
\[v=u+iB^{-1}u_t \]
where $B^2=Id-\Delta.$  We can rewrite \eqref{kgh} in the form
\begin{eqnarray}\label{bit}
\begin{cases}
iv_t-Bv-B^{-1} \left( [V \ast | \text{Re} \ v|^2 ] \text{Re} \ v \right) =0\\
v(0)=u(0) + i B^{-1} \left( u_t(0) \right) = F_N + F^N,
\end{cases}
\end{eqnarray} where  
\[F_N:= f_N+i B^{-1}(g_N)\in H^1(\R^3) \quad \text{and} \quad F^{N}:=f^{N}+i B^{-1}(g^N)\in M^{p_\gamma, p'_\gamma}(\R^3). \]
 By \eqref{tbu} and Proposition \ref{fp} (with $\sigma=s=1$), we  have 
\begin{eqnarray}\label{li}
\|F_N\|_{H^1}\lesssim N^{\frac{\theta}{1-\theta}}, \quad   \|F^N\|_{M_1^{p_\gamma, p_\gamma'}} \lesssim \frac{1}{N},
\end{eqnarray}
and 
more generally,
\begin{eqnarray}\label{mg}
\| \left( u(t, \cdot), \partial_t u(t, \cdot) \right)\|_{H^1\times L^2+M_1^{p_\gamma,p'_\gamma}\times M^{p_\gamma,  p'_\gamma}} \approx \|v(t, \cdot)\|_{H^1 + M^{p_\gamma, p'_\gamma}_1}. 
\end{eqnarray}
The Hamiltonian of \eqref{bit}  (cf.  \cite[Equation IV.2.7]{Bourgain1999})  is formally given by the formula
\begin{eqnarray}\label{hamil}
H(v)= \int_{\mathbb R^3} \left[ \frac{1}{2} |Bv|^2 + \frac{1}{4} (V \ast |\text {Re}  \ v|^2) |\text{Re} \ v|^2 \right] dx.
\end{eqnarray}
However,  we cannot use it to control the  full solution $v,$ since \eqref{hamil} does not make sense for general $v(t, \cdot) \in H^1(\R^3) + M_1^{p_\gamma, p'_\gamma}(\R^3).$
By Proposition \ref{bm} the linear evolution $e^{-it B}F^N$ never blows up in modulation spaces on any bounded time interval $I=[0,T].$  Specifically,  by  Proposition \ref{bm} we have 
\begin{eqnarray*}
\|e^{-it B}\|_{M^{p_\gamma, p'_\gamma}_1 \to M_1^{p_\gamma, p'_\gamma}} \lesssim 1 \quad \text{for all} \ t \in [0,T].
\end{eqnarray*}
Therefore,  we are left to control
\begin{eqnarray}\label{upe}
\tilde{v}:=v- e^{-it B} F^N
\end{eqnarray}
in  $H^1$ on $I=[0,T]$ (see \cite[Equation IV 2.21]{Bourgain1999}) using the Hamiltonian.  
By \eqref{upe} and  \eqref{bit},  we have $\tilde{v}(0)=v(0)-F^N=F_N.$ Define $$I(t)= H(\tilde{v}(t)).$$
Given $\gamma \in (2,3),$ we can choose $r_1, p_1 \in (1,3]$ such that
\begin{equation} \label{ufu}
\frac{1}{p_1}+\frac{1}{q_1}=1, \quad \quad  \frac{1}{r_1}+\frac{\gamma}{3}-1= \frac{1}{q_1} \quad  \text{and} \quad r_1<q_1.
\end{equation}
By \eqref{ufu} and    H\"older inequality,  Propositions \ref{hls} and \ref{up},  and \eqref{li},   we obtain 
\begin{align}\label{s1}
& I(0)   =  H(F_N)= \int_{\mathbb R^3} \left[ \frac{1}{2} |B F_N|^2 + \frac{1}{4} (V \ast |\text {Re}   F_N|^2) |\text{Re} F_N|^2 \right] dx \nonumber \\
 & \lesssim  \|F_N\|^2_{H^1} + \|  |\cdot |^{-\gamma} \ast |\text{Re}  F_N|^2\|_{L^{q_1}} \||\text {Re} F_N|^2\|_{L^{p_1}} \nonumber \\ 
& \lesssim \|F_N\|^2_{H^1} +   \| |F_N|^2\|_{L^{r_1}} \|F_N\|^2_{L^{2p_1}} \nonumber\\
& \lesssim \|F_N\|^2_{H^1} +  \|F_N\|_{H^1}^4 \nonumber\\
& \lesssim N^{\frac{2\theta}{1-\theta}} + N^{\frac{4\theta}{1-\theta}} \nonumber\\
& \lesssim  N^{\frac{4\theta}{1-\theta}}.
\end{align}
We aim to  estimate the time $T$ (as a function of $N$) that preserves \eqref{s1}, i.e. $T$  such that for all $0\leq t \leq T$ we have
\begin{eqnarray}\label{tbt}
I(t)& = & \int_{\R^3} \left[ \frac{1}{2} |B\tilde{v}(t,x)|^2 + \frac{1}{4} (V \ast |\text {Re}  \ \tilde{v}(t,x)|^2) |\text{Re} \ \tilde{v}(t,x)|^2 \right] dx \nonumber \\
&  \leq  & 2I(0) \lesssim N^{\frac{4\theta}{1-\theta}}.
\end{eqnarray} 
Exploiting ideas from  \cite[Appendix B]{Pattakos2021} and \cite[Equation IV 2.26]{Bourgain1999},  we may  formally obtain
\begin{eqnarray*}
I'(t)= \text{Im} \left\langle  B \tilde{v},  (V \ast |\text {Re} \ v|^2)\text{Re} \ v- (V \ast |\text {Re} \ \tilde{v}|^2)\text{Re} \ \tilde{v} \right\rangle. 
\end{eqnarray*}
By Cauchy-Schwartz inequality and \eqref{tbt} we have
\begin{eqnarray}\label{ss}
|I'(t)|  & \leq & \|B \tilde{v}\|_{L^2} \| (V \ast |\text {Re} \ v|^2)\text{Re} \ v- (V \ast |\text {Re} \ \tilde{v}|^2)\text{Re} \ \tilde{v} \|_{L^2} \nonumber\\
& \lesssim & I(t)^{\frac{1}{2}}  \| (V \ast |\text {Re} \ v|^2)\text{Re} \ v- (V \ast |\text {Re} \ \tilde{v}|^2)\text{Re} \ \tilde{v} \|_{L^2} \nonumber \\
& \lesssim & N^{\frac{2\theta}{1-\theta}}  \| (V \ast |\text {Re} \ v|^2)\text{Re} \ v- (V \ast |\text {Re} \ \tilde{v}|^2)\text{Re} \ \tilde{v} \|_{L^2}.
\end{eqnarray}
In view of \eqref{pehn} and \eqref{upe}, we  have 
\begin{eqnarray} \label{a}
&& \left|(V\ast | \text{Re} \ v|^{2}) \text{Re} \ v- (V\ast | \text{Re} \ \tilde{v}|^{2}) \text{Re} \ \tilde{v}\right|  \nonumber\\
 & \leq  &  \left|   (V \ast (| \text{Re} \ v|^{2})e^{-itB} F^N \right|  + \left|(V \ast (| \text{Re} \ v|^{2}- | \text{Re} \ \tilde{v} \  |^{2})) \text{Re} \ \tilde{v}\right|.
\end{eqnarray}
In view of \eqref{uie} and similar argument as in the proof of Theorem \ref{mt1}
and Proposition \ref{fp} (with $\sigma=s=1$) and Proposition \ref{up} \eqref{eml},  and \eqref{li}, we obtain
\begin{eqnarray}\label{b}
\left\| (V \ast (| \text{Re} \ v|^{2})e^{-itB} F^N \right\|_{L^2}  & \lesssim &  \|e^{-itB}F^N \|_{L^{\frac{18}{9-2\gamma}}} \|(V\ast | \text{Re} v|^2) \|_{L^{\frac{9}{\gamma}}} \\
& \lesssim & \|e^{-itB}F^N\|_{M^{p_\gamma, p'_\gamma}} \| \text{Re} \ v\|^2_{L^{p_\gamma}} \nonumber \\
& \lesssim &  \|F^N\|_{M_1^{p_\gamma, p'_\gamma}} \|e^{-itB} F^N\|_{L^{p_\gamma}}^{2} \nonumber \\
& \lesssim &  \|F^N\|_{M_1^{p_\gamma, p_\gamma}} \|e^{-itB} F^N\|_{M^{p_\gamma, p'_\gamma}}^{2} \nonumber\\
& \lesssim & N^{-3}.
\end{eqnarray}
By H\"older inequality,  Proposition \ref{hls} and Proposition \ref{up}\eqref{up1}, we have 
\begin{eqnarray*}
 \|(V \ast (| \text{Re} \ v|^{2}- | \text{Re} \ \tilde{v} \  |^{2})) \text{Re} \ \tilde{v}\|_{L^2}  & \leq  & \|(V \ast (| \text{Re} \ v|^{2}- | \text{Re} \ \tilde{v} \  |^{2})) \|_{L^{\frac{9}{\gamma}}} \| \text{Re} \ \tilde{v}\|_{L^{\frac{18}{9-2\gamma}}} \\
 & \lesssim &  \| | \text{Re} \ v|^{2}- | \text{Re} \ \tilde{v} \  |^{2} \|_{L^{\frac{9}{9-2\gamma}}} \|  \ \tilde{v}\|_{H^{1}}.
\end{eqnarray*}
In view of \eqref{upe},  we obtain
\begin{eqnarray*}
\left| | \text{Re} \ v|^{2}- | \text{Re} \ \tilde{v} \  |^{2} \right| & = & \left|  \text{Re}  \ v  \left(  \overline{ \text{Re}  \ v - \text{Re} \  \tilde{v}} \right) + \overline{\text{Re}  \ \tilde{v}} \left( \text{Re}  \ v - \text{Re}  \ \tilde{v} \right) \right|\\
& = & \left|  \text{Re}  \ v \   \overline{ e^{-itB} F^N } + \overline{\text{Re}  \ \tilde{v}}  \  e^{-itB} F^N  \right|\\
& \lesssim &  |e^{-itB}F^N|^2 +  |\tilde{v}| |e^{-itB}F^N|.
\end{eqnarray*} 
Using the above two inequalities,  Proposition \ref{up}\eqref{eml} and \eqref{li},   for large $N$,  we have 
\begin{eqnarray}\label{c}
&&  \|(V \ast (| \text{Re} \ v|^{2}- | \text{Re} \ \tilde{v} \  |^{2})) \text{Re} \ \tilde{v}\|_{L^2} \nonumber \\
 & \lesssim &  \left( \|e^{-itB} F^N\|^2_{L^{\frac{18}{9-2\gamma}}} +  \|e^{-itB} F^N\|_{L^{\frac{18}{9-2\gamma}}} \|\tilde{v}\|_{L^{\frac{18}{9-2\gamma}}}\right) \| \tilde{v} \|_{H^1}  \nonumber\\
 & \lesssim &  \left( \| e^{-itB} F^N\|^2_{M^{\frac{18}{9-2\gamma},  \left(\frac{18}{9-2\gamma} \right)'} } + \| e^{-itB} F^N\|_{M^{\frac{18}{9-2\gamma},  \left(\frac{18}{9-2\gamma} \right)'} } \| \tilde{v} \|_{H^1} \right)\| \tilde{v} \|_{H^1} \nonumber \\
 & \lesssim & \left( \| F^N\|^2_{M_1^{\frac{18}{9-2\gamma},  \left(\frac{18}{9-2\gamma} \right)'} } + \| F^N\|_{M_1^{\frac{18}{9-2\gamma},  \left(\frac{18}{9-2\gamma} \right)'} } I(t)^{\frac{1}{2}} \right)I(t)^{\frac{1}{2}}  \nonumber \\
  & \lesssim & \left( N^{-2} + N^{-1} N^{\frac{2\theta}{1-\theta}} \right) N^{\frac{2\theta}{1-\theta}}  \nonumber \\
  & \lesssim & N^{-1+\frac{4\theta}{1-\theta}}.
\end{eqnarray}
Combining \eqref{a}, \eqref{b} and \eqref{c},  we have 
\begin{eqnarray*}
\|(V\ast | \text{Re} \ v|^{2}) \text{Re} \ v- (V\ast | \text{Re} \ \tilde{v}|^{2}) \text{Re} \ \tilde{v}\|_{L^2} \lesssim N^{-1+\frac{4\theta}{1-\theta}} + N^{-3}.
\end{eqnarray*}
Using  the mean value theorem and the above estimate in \eqref{ss}, for $0\leq t \leq T$ and some $\tau \in [0, t],$ we have 
\begin{eqnarray*}
\left| I(t)-I(0) \right| & \leq & T \left| I'(\tau) \right| \\\
& \lesssim & T N^{\frac{2\theta}{1-\theta}} \left( N^{-1+\frac{4\theta}{1-\theta}} + N^{-3} \right) \nonumber \\
& \lesssim & T\left( \frac{1}{N^{1- \frac{6\theta}{1-\theta}}} + \frac{1}{N^{3-\frac{2\theta}{1-\theta}}} \right). \nonumber
\end{eqnarray*}
Thus,  we have 
\begin{eqnarray}\label{tbtd}
I(t) \lesssim  I(0) + T\left( \frac{1}{N^{1- \frac{6\theta}{1-\theta}}} + \frac{1}{N^{3-\frac{2\theta}{1-\theta}}} \right)
\end{eqnarray}
For \eqref{tbt}  to be true it suffices that the last  term of \eqref{tbtd}
satisfies  $$T\left( \frac{1}{N^{1- \frac{6\theta}{1-\theta}}} + \frac{1}{N^{3-\frac{2\theta}{1-\theta}}} \right) \lesssim I(0).$$
Taking \eqref{s1} into account,  we  impose the following condition
$$T\left( \frac{1}{N^{1- \frac{6\theta}{1-\theta}}} + \frac{1}{N^{3-\frac{2\theta}{1-\theta}}} \right) \lesssim  N^{\frac{4\theta}{1-\theta}}.$$
In other words,  for \eqref{tbt}  to be true  we must have the following condition
\begin{eqnarray}
T\left( \frac{1}{N^{1- \frac{2\theta}{1-\theta}}} + \frac{1}{N^{3+\frac{2\theta}{1-\theta}}} \right) \lesssim 1.
\end{eqnarray}
Since $N$ is going to be large  we want both the exponents of $N$ in the last expression are positive.  Therefore,  we require that
\begin{eqnarray}\label{rc}
 \frac{2\theta}{1-\theta} <1.
\end{eqnarray}
Since  
$$ \frac{27-2\gamma}{54}<\frac{1-\theta}{2}+ \frac{\theta}{p_\gamma}< \frac{1}{2},$$
we get $\theta< \frac{1}{3}$ and so \eqref{rc} is satisfied.  Also,  we have 
\[ 1- \frac{2\theta}{1-\theta} \leq  3+ \frac{2\theta}{1- \theta}.\]
Hence,  \eqref{tbt} holds for $$T \sim N^{1- \frac{2\theta}{1-\theta}},   \quad \text{i.e.} \ N=N(T) \sim T^{\left(1- \frac{2\theta}{1-\theta} \right)^{-1}}.$$  Thus,  in view of  \eqref{bit},  \eqref{upe} and \eqref{li},  for $0\leq t \leq T$, we obtain
\begin{eqnarray*}
\| v(t)- e^{-it B} v(0)\|_{H^1}  & = & \| \tilde{v}(t)- e^{-it B} F_N\|_{H^1} \leq \| \tilde{v}(t)\|_{H^1} + \| e^{-it B} F_N\|_{H^1} \\
& \lesssim &  I(t)^{\frac{1}{2}} + \|F_N\|_{H^1}  \lesssim  N^{\frac{2\theta}{1-\theta}} + N^{\frac{\theta}{1-\theta}}\\
& \lesssim  & N^{\frac{2\theta}{1-\theta}} \sim T^{\frac{2\theta}{1-3\theta}}
\end{eqnarray*}
Consequently,
 \begin{eqnarray*}
 \|\left( u(t), u_t(t) \right) - \mathbb K(t) \left( u(0),  u_t(0) \right) \|_{H^1 \times L^2} \lesssim  (1+t)^{\frac{2\theta}{1-3\theta}}.
 \end{eqnarray*}
 This completes the proof of Theorem \ref{mt1}.
\end{proof}

\bibliographystyle{amsplain}
\bibliography{kgh}
\end{document}